\documentclass[11pt,a4paper]{article}
\ifdefined\pdfminorversion\pdfminorversion=4\fi
\ifdefined\pdfobjcompresslevel\pdfobjcompresslevel=0\fi
\IfFileExists{lmodern.sty}{\usepackage[T1]{fontenc}\usepackage{lmodern}}{}
\IfFileExists{microtype.sty}{\usepackage[nopatch=eqnum]{microtype}}{}
\usepackage{amsmath,amssymb} 
\usepackage{amsthm}          
\usepackage{cite}            
\usepackage{geometry}
\usepackage{indentfirst}    
\usepackage{xcolor}

\newtheorem{theorem}{Theorem}[section]
\newtheorem{lemma}[theorem]{Lemma}
\newtheorem{claim}[theorem]{Claim}

\theoremstyle{remark}
\newtheorem{remark}[theorem]{Remark}
\theoremstyle{plain}

\title{Edge-Disjoint Cliques above the Tur\'an Number}

\author{Hanzhi Bai\thanks{School of Mathematics, Shandong University, Jinan 250100, China. Email: \mbox{bhz@mail.sdu.edu.cn}.}
\and Jialin He\thanks{School of Mathematical Sciences, Key Laboratory of MEA, and Shanghai Key Laboratory of PMMP, East China Normal University, Shanghai 200241, China. Email: \mbox{jlhe@math.ecnu.edu.cn}.}
\and Fuhong Ma\thanks{School of Mathematics, Shandong University of Technology, Zibo 255000, China.\newline Email: \mbox{mafuhongsdnu@163.com}.}
\and Jin Yan\thanks{School of Mathematics, Shandong University, Jinan 250100, China. Email: \mbox{yanj@sdu.edu.cn}.}}
\date{}

\usepackage[colorlinks=true,linkcolor=red,citecolor=blue,urlcolor=blue]{hyperref}
\hypersetup{pdftitle={Edge-Disjoint Cliques above the Tur\'an Number},pdfauthor={Hanzhi Bai, Jialin He, Fuhong Ma, Jin Yan}}

\let\OriginalCitepunct\citepunct
\let\OriginalCitedash\citedash
\renewcommand{\citepunct}{\textcolor{blue}{\OriginalCitepunct}}
\renewcommand{\citedash}{\textcolor{blue}{\OriginalCitedash}}
\mathchardef\citeprepenalty=10000
\makeatletter
\DeclareRobustCommand{\eqref}[1]{\textcolor{red}{\textup{\tagform@{\ref{#1}}}}}
\makeatother

\begin{document}
\maketitle

\begin{abstract}
Fix $r\ge3$, and let $t_{r-1}(n)$ be the Tur\'an number for $K_r$. Let $\phi_r(n,k)$ be the minimum, over all $n$-vertex graphs with $t_{r-1}(n)+k$ edges, of the maximum number of edge-disjoint copies of $K_r$. We answer a question of Balogh and Wigal by proving that, for positive integers $k=k(n)$ with $t_{r-1}(n)+k\le\binom n2$, $\phi_r(n,k)=(1-o(1))k$ if and only if $k=o(n^2)$. When $k/n^2$ tends to a positive constant, we give an explicit asymptotic upper bound below one for $\phi_r(n,k)/k$. Gy\H{o}ri proved that $k-\phi_r(n,k)=O_r(k^2/n^2)$ when $k=o(n^2)$. For $k/n\to\infty$ and $k/n^2\to0$, we sharpen this bound to
\[
 k-\phi_r(n,k)=(A_r+o(1))\frac{k^2}{n^2},
\]
where $A_r$ is an explicit constant and $A_3=4/9$. The same coefficient holds for fractional packings. The matching construction starts with a clique joined to a complete $(r-1)$-partite graph with unequal part sizes. For the lower bound on the packing number, we first prove a uniform estimate for fractional packings. We then obtain the integer bound by reducing the gap between fractional and integer packings to that of a graph on $O_r(k/n)$ vertices.
\end{abstract}


\section{Introduction}

Let $T_{r-1}(n)$ be the balanced complete $(r-1)$-partite graph on $n$ vertices, and let $t_{r-1}(n)$ be its number of edges. By Tur\'an's theorem, every $n$-vertex graph with more than $t_{r-1}(n)$ edges contains a copy of $K_r$. How many edge-disjoint copies must such a graph contain if it has $k$ edges above the threshold? Write $e(G)$ for the number of edges in $G$. Let $\nu_r(G)$ be the maximum number of pairwise edge-disjoint copies of $K_r$ in $G$, and let
\[
 \phi_r(n,k)=\min\{\nu_r(G): |V(G)|=n,\ e(G)=t_{r-1}(n)+k\}.
\]
Throughout, $r\ge3$ is fixed. We call $k=e(G)-t_{r-1}(n)$ the excess of an $n$-vertex graph $G$. An integer $k$ is feasible if $0\le k\le\binom n2-t_{r-1}(n)$. Adding $k$ edges to $T_{r-1}(n)$ gives $\phi_r(n,k)\le k$, since every $K_r$ in the resulting graph uses an added edge. We study the loss $k-\phi_r(n,k)$ from this upper bound.

For $r=3$, Erd\H{o}s~\cite[p.~98]{Erdos1971} and Sauer~\cite{Sauer} studied how long the equality $\phi_3(n,k)=k$ persists as $k$ increases. For sufficiently large $n$, Gy\H{o}ri determined the exact initial range, which ends at $2n-10$ for odd $n$~\cite[Theorem~9 and Remark~2]{Gyori1992} and at $3n/2-5$ for even $n$~\cite{Gyori1988}. For fixed $r\ge4$, he proved $\phi_r(n,k)=k$ for $0\le k\le3\lfloor(n+1)/(r-1)\rfloor-5$ and sufficiently large $n$~\cite{Gyori1991}. More generally, for any positive feasible $k=o(n^2)$, Gy\H{o}ri proved that
\begin{equation}\label{eq:gyori}
 k-\phi_r(n,k)=O_r(k^2/n^2).
\end{equation}
See~\cite[Theorem~8]{Gyori1992} and~\cite[Theorem~2.2]{LPS}.

Balogh and Wigal~\cite{BaloghWigal2025} proved $\phi_r(n,k)\ge(2/r-o(1))k$ over the whole feasible range. For triangles, they gave examples in which $\phi_3(n,k)/k$ stays bounded away from one while $k/n^2$ tends to a positive constant. They asked for which positive feasible integers $k=k(n)$ the ratio $\phi_r(n,k)/k$ tends to one. Our first result answers this question and gives an explicit upper bound whenever $k/n^2$ has a positive limit.

\begin{theorem}\label{thm:range}
Fix $r\ge3$, and let $k=k(n)$ be a positive feasible integer. Then $\phi_r(n,k)=(1-o(1))k$ if and only if $k=o(n^2)$.
Moreover, if $k/n^2\to\delta$ with $0<\delta\le1/[2(r-1)]$, then
\[
 \frac{\phi_r(n,k)}k
 \le1-\frac{r-2}{r}\,\frac{1-\sqrt{1-2(r-1)\delta}}{1+\sqrt{1-2(r-1)\delta}}+o(1).
\]
\end{theorem}

The constant in this upper bound is strictly less than one. At $\delta=1/[2(r-1)]$, it equals $2/r$ and matches the general lower bound.
Theorem~\ref{thm:range} determines exactly when $k-\phi_r(n,k)=o(k)$. We next determine the leading term of this loss when $k/n\to\infty$ and $k/n^2\to0$.

A fractional $K_r$-packing of $G$ consists of weights $w_K\ge0$, one for each copy $K$ of $K_r$ in $G$, satisfying $\sum_{K:\,e\in E(K)}w_K\le1$ for every $e\in E(G)$.
The sum on the left is the load on $e$. Let $\nu_r^*(G)=\max\sum_Kw_K$, where the maximum is over all such packings, and let
\[
 \phi_r^*(n,k)=\min\{\nu_r^*(G): |V(G)|=n,\ e(G)=t_{r-1}(n)+k\}.
\]
Write $n\ll k\ll n^2$ for $k/n\to\infty$ and $k/n^2\to0$. Define
\[
 A_r=(r-1)^2\left(\sum_{j=1}^r\frac1{j^2}-\frac32+\frac1r-\frac1{2r(r-1)}\right).
\]

\begin{theorem}\label{thm:second-order}
Fix $r\ge3$. If $k=k(n)$ is feasible and $n\ll k\ll n^2$, then
\[
 \phi_r(n,k)=k-(A_r+o(1))\frac{k^2}{n^2}.
\]
The same asymptotic formula holds for fractional packings:
\[
 \phi_r^*(n,k)=k-(A_r+o(1))\frac{k^2}{n^2}.
\]
\end{theorem}

In particular, $A_3=4/9$, $A_4=19/16$, and $A_5=499/225$. The formula also gives $A_r/r^2\to\pi^2/6-3/2$ as $r\to\infty$. The error terms in Theorem~\ref{thm:second-order} are for fixed $r$. We obtain the fractional lower bound from the following estimate, which holds over a wider range.

\begin{theorem}\label{thm:uniform}
For every fixed $r\ge3$ there are constants $c_r,B_r>0$ such that every graph $G$ on $n\ge1$ vertices with $e(G)=t_{r-1}(n)+k$ and $0\le k\le c_r n^2$ satisfies
\[
 \nu_r^*(G)\ge k-A_r\frac{k^2}{n^2}-B_r\frac{k^{5/2}}{n^3}-B_r.
\]
\end{theorem}

An additional clique restriction can give the full lower bound $k$. Gy\H{o}ri and Keszegh~\cite{GyoriKeszegh2017} proved that a $K_4$-free graph with $t_2(n)+k$ edges contains $k$ edge-disjoint triangles. Liu and Zeng~\cite{LiuZeng} extended this to $K_{r+1}$-free graphs for each fixed $r$ and all sufficiently large $n$. Thus, for sufficiently large $n$, every $n$-vertex graph $G$ with $e(G)=t_{r-1}(n)+k$ and $\nu_r(G)<k$ contains a $K_{r+1}$.

Related work concerns clique coverings and decompositions~\cite{EGP,Bollobas,GT}. For triangles, Kr\'al' et~al.~\cite{KLMP} proved \mbox{$\phi_3(n,k)\ge2k/3-o(n^2)$}. Blumenthal et~al.~\cite{BLPPPV} improved this to \mbox{$\phi_3(n,k)\ge2k/3-O(1)$}. For decompositions allowing cliques of different orders with cost $|V(C)|-1$ per clique $C$, see~\cite{BHKNW,Zhang}.

To study fractional packings, we also use their dual covering problem. A $K_r$-cover of $G$ is a set of edges that meets every copy of $K_r$ in $G$. A fractional $K_r$-cover assigns weights $w_{uv}\in[0,1]$ to the edges $uv\in E(G)$ such that $\sum_{uv\in E(K)}w_{uv}\ge1$ for every copy $K$ of $K_r$ in $G$.
Its cost is $\sum_{uv\in E(G)}w_{uv}$. By linear programming duality, $\nu_r^*(G)=\min_w\sum_{uv\in E(G)}w_{uv}$, where the minimum is over all fractional $K_r$-covers of $G$. Thus any fractional cover of cost $W$ gives \mbox{$\nu_r(G)\le\nu_r^*(G)\le W$}.

All graphs are finite and simple, and constants with subscript $r$ depend only on $r$. Unless otherwise stated, all asymptotic statements are as $n\to\infty$.

In Section~\ref{sec:balanced}, we prove Theorem~\ref{thm:range} and give the construction for the upper bounds in Theorem~\ref{thm:second-order}. Section~\ref{sec:fractional} establishes Theorem~\ref{thm:uniform} and the fractional assertion of Theorem~\ref{thm:second-order}. Finally, we prove the integer assertion of Theorem~\ref{thm:second-order} in Section~\ref{sec:integer}.

\section{Proof of Theorem~\ref{thm:range} and the Upper Bound in Theorem~\ref{thm:second-order}}\label{sec:balanced}

We use one family of graphs, obtained by joining a clique to a complete $(r-1)$-partite graph. For $r=3$, a construction with balanced parts already appears in Gy\H{o}ri~\cite[Example~1]{Gyori1992}. Balanced part sizes and equal weights on edges from the clique to the parts give Theorem~\ref{thm:range}. We then choose unequal part sizes and weights to obtain the upper bounds in Theorem~\ref{thm:second-order}. In both cases, we compare the excess with the cost of a fractional cover.

A graph in this family has a partition into a clique $S$ of size $t$ and independent sets $A_1,\ldots,A_{r-1}$ of sizes $b_1,\ldots,b_{r-1}$, with $t+\sum_i b_i=n$ and all edges between distinct sets. Empty parts are allowed. To obtain a fractional cover, we give edges inside $S$ weight $2/[r(r-1)]$, edges between distinct $A_i$ weight zero, and edges from $S$ to $A_i$ weight $x_i$, where $0\le x_1\le\cdots\le x_{r-1}$ and $\sum_i x_i=1$. Every copy of $K_r$ meets $S$. If it has $j$ vertices in $S$, its other $r-j$ vertices lie in distinct parts. The ordering of the $x_i$ therefore gives total weight at least
\[
 \frac{j(j-1)}{r(r-1)}+j\sum_{i=1}^{r-j}x_i.
\]
Thus these weights form a fractional cover whenever this expression is at least one for $1\le j\le r$. We use the convention that an empty sum is zero.

The edge count uses the bounded rounding term in the Tur\'an number. For an integer $m\ge0$, write $m=(r-1)\ell+d$ with $0\le d<r-1$, and let
\[
 \eta_{r-1}(m)=\frac{r-2}{2(r-1)}m^2-t_{r-1}(m)
 =\frac{d(r-1-d)}{2(r-1)}\in[0,(r-1)/8].
\]
The excess $k'$ of this graph is
\[
 k'=\binom n2-\sum_i\binom{b_i}{2}-t_{r-1}(n)
 =\frac{n^2}{2(r-1)}-\frac12\sum_i b_i^2-\frac t2+\eta_{r-1}(n).
\]
The cost of the fractional cover is
\[
 W=\frac{2\binom t2}{r(r-1)}+t\sum_i b_ix_i.
\]
Since $\sum_i(b_i+tx_i)=n$, completing the square gives
\begin{equation}\label{eq:sharp-construction}
\begin{aligned}
 k'-W
 &=\left(\frac12\sum_i x_i^2-\frac1{r(r-1)}\right)t^2
   -\left(\frac12-\frac1{r(r-1)}\right)t+\eta_{r-1}(n)\\
 &-\frac12\sum_{i=1}^{r-1}\left(b_i-\frac n{r-1}+tx_i\right)^2.
\end{aligned}
\end{equation}
The graph can be adjusted to any prescribed feasible excess $k\ge k'$ by adding $k-k'$ missing edges of weight one. Every new copy of $K_r$ contains an added edge, so the resulting fractional cover has cost $W+k-k'$. The difference between the new excess and this cost remains $k-(W+k-k')=k'-W$. Thus
\begin{equation}\label{eq:construction-adjust}
 \phi_r(n,k)\le\phi_r^*(n,k)\le k-(k'-W).
\end{equation}

\smallskip
\begin{proof}[Proof of Theorem~\ref{thm:range}]
For the first-order bound, we use the balanced member with $t=a$ and $x_i=1/(r-1)$ for every $i$, so the part sizes $b_i$ differ by at most one. Its excess and cover cost are denoted by $k_a$ and $W_a$. The covering condition holds because
\[
 \frac{j(j-1)}{r(r-1)}+\frac{j(r-j)}{r-1}
 =1+\frac{(j-1)(r-j)}r\ge1.
\]
The edge count and~\eqref{eq:sharp-construction} give
\begin{gather*}
 k_a=\binom a2+a(n-a)+t_{r-1}(n-a)-t_{r-1}(n)
 =\frac{a(2n-a)}{2(r-1)}-\frac a2+\eta_{r-1}(n)-\eta_{r-1}(n-a),\\
 k_a-W_a=\frac{r-2}{2r(r-1)}a^2-\left(\frac12-\frac1{r(r-1)}\right)a
 +\eta_{r-1}(n)-\eta_{r-1}(n-a).
\end{gather*}
The bounded rounding terms give $k_a-W_a=\frac{r-2}{2r(r-1)}a^2+O_r(a+1)$.
Moving a vertex of a largest part into $S$ leaves the $A_i$ balanced and adds the edges to the other vertices of that part. Hence $k_{a+1}-k_a=\lceil(n-a)/(r-1)\rceil-1$ for $a<n$. The sequence $k_a$ is nondecreasing from $k_0=0$ to $k_n=\binom n2-t_{r-1}(n)$, with increments at most $n/(r-1)$.

Choose the largest $a$ with $k_a\le k$. If $a<n$, then $k<k_{a+1}$; if $a=n$, then $k=k_n$. In either case, $0\le k-k_a\le n/(r-1)$. If $k/n^2\to\delta>0$, feasibility gives $\delta\le1/[2(r-1)]$. The formula for $k_a$ yields $(a/n)(2-a/n)\to2(r-1)\delta$. Since $0\le a/n\le1$, we obtain $a/n\to1-\sqrt{1-2(r-1)\delta}$. By~\eqref{eq:construction-adjust},
\begin{align*}
 \frac{\phi_r(n,k)}k
 &\le\frac{\phi_r^*(n,k)}k
 \le1-\frac{(r-2)(1-\sqrt{1-2(r-1)\delta})^2}{2r(r-1)\delta}+o(1)\\
 &=1-\frac{r-2}{r}\,\frac{1-\sqrt{1-2(r-1)\delta}}{1+\sqrt{1-2(r-1)\delta}}+o(1).
\end{align*}
The constant on the right is strictly less than one.

If $k=k(n)$ is positive and feasible but is not $o(n^2)$, there is a subsequence of values of $n$ along which $k/n^2$ tends to a positive limit. The preceding bound rules out $\phi_r(n,k)/k\to1$. Conversely, if $k=o(n^2)$, then~\eqref{eq:gyori} gives $\phi_r(n,k)/k\to1$. Taking $a=0$ in~\eqref{eq:construction-adjust} gives $\phi_r(n,k)\le\phi_r^*(n,k)\le k$, so $\phi_r^*(n,k)/k\to1$ as well.
\end{proof}

\begin{remark}
At $\delta=1/[2(r-1)]$, the upper bound is $2/r$ and matches the general lower bound. For complete graphs, the same limit follows from Wilson's theorem~\cite{Wilson1976}. For intermediate densities, counting the edges in a packing also gives $\nu_r(G)\le e(G)/\binom r2$. Thus, if $k/n^2\to\delta>0$, then
\[
 \frac{\phi_r(n,k)}{k}\le\frac{(r-2)/(2(r-1))+\delta}{\binom r2\delta}+o(1).
\]
This bound may be smaller than that in Theorem~\ref{thm:range}, so one may take the minimum of the two bounds. Determining the sharp asymptotic value at intermediate densities remains open in general.
\end{remark}

\smallskip
\begin{proof}[Proof of Theorem~\ref{thm:second-order} (Upper Bound)]
We use the same family with $n\ll k\ll n^2$ and \mbox{$t=\lfloor(r-1)k/n\rfloor$}, so $t\to\infty$ and $t=o(n)$. The equal weights $x_i=1/(r-1)$ make the covering bound strict for $1<j<r$. We therefore choose the $x_i$ so that this bound equals one for every $j$. With $j=r-\ell$, this requires, for each $0\le\ell\le r-1$,
\begin{equation}\label{eq:prefix-weights}
 \sum_{i=1}^{\ell}x_i=\frac1{r-\ell}-\frac{r-\ell-1}{r(r-1)}.
\end{equation}
Subtracting consecutive partial sums gives, for each $1\le i\le r-1$,
\begin{equation}\label{eq:optimal-weights}
 x_i=\frac1{(r-i)(r-i+1)}+\frac1{r(r-1)}.
\end{equation}
These weights are positive and increasing. By~\eqref{eq:prefix-weights}, they sum to one and satisfy every covering condition.

With these weights fixed, the last term of~\eqref{eq:sharp-construction} is largest when $b_i=n/(r-1)-tx_i$, since all its squares then vanish. These real sizes sum to $n-t$ and are positive for large $n$. Round them down and distribute the remaining vertices among parts whose prescribed sizes are nonintegral, adding at most one vertex to each. Each size changes by less than one, so rounding lowers $k'-W$ by at most $(r-1)/2$.

We check that the rounded graph has at most the prescribed number of edges. For a fixed $t$, the balanced member of the family has the most edges, since it minimizes $\sum_i\binom{b_i}{2}$. Its excess is at most $tn/(r-1)$ by the increment bound for $k_a$. Hence $k'\le tn/(r-1)\le k$, so~\eqref{eq:construction-adjust} applies to the rounded graph.

To evaluate the coefficient of $t^2$, use $1/[j(j+1)]=1/j-1/(j+1)$ and $\sum_{j=1}^{r-1}1/[j(j+1)]=(r-1)/r$. Then
\begin{align*}
 \sum_{j=1}^{r-1}\frac1{j^2(j+1)^2}
 &=\sum_{j=1}^{r-1}\left(\frac1{j^2}+\frac1{(j+1)^2}-\frac2{j(j+1)}\right)
 =\left(\sum_{j=1}^r\frac1{j^2}-\frac1{r^2}\right)
 +\left(\sum_{j=1}^r\frac1{j^2}-1\right)-\frac{2(r-1)}r\\
 &=2\sum_{j=1}^r\frac1{j^2}-3+\frac2r-\frac1{r^2}.
\end{align*}
Substituting $j=r-i$ in~\eqref{eq:optimal-weights} gives
\begin{align*}
 \sum_{i=1}^{r-1}x_i^2
 &=\sum_{j=1}^{r-1}\left(\frac1{j(j+1)}+\frac1{r(r-1)}\right)^2
 =\sum_{j=1}^{r-1}\frac1{j^2(j+1)^2}
 +\frac2{r(r-1)}\frac{r-1}r+\frac1{r^2(r-1)}\\
 &=2\sum_{j=1}^r\frac1{j^2}-3+\frac2r+\frac1{r(r-1)}.
\end{align*}
Consequently, $\tfrac12\sum_i x_i^2-1/[r(r-1)]=A_r/(r-1)^2$, which is positive because $\sum_i x_i^2\ge1/(r-1)>2/[r(r-1)]$.

By~\eqref{eq:sharp-construction}, the rounding bound, and $t=(r-1)k/n+O(1)$,
\[
 k'-W=\frac{A_r}{(r-1)^2}t^2+O_r(t+1)
 =A_r\frac{k^2}{n^2}+O_r(k/n+1).
\]
Since $k/n\to\infty$, equation~\eqref{eq:construction-adjust} yields
\begingroup
\makeatletter
\def\equation@qed{%
  \gdef\alt@tag{%
    \global\let\alt@tag\@empty
    \hbox{\tagform@\theequation\quad\qedsymbol}%
    \setbox\z@
  }%
}
\makeatother
\begin{equation}\label{eq:sharp-upper}
 \phi_r(n,k)\le\phi_r^*(n,k)
 \le k-(A_r-o(1))\frac{k^2}{n^2}.\qedhere
\end{equation}
\endgroup
\end{proof}

\begin{remark}
For $r=3$, the weights are $x_1=1/3$ and $x_2=2/3$. Writing the part sizes as $(n-t)/2+d$ and $(n-t)/2-d$, identity~\eqref{eq:sharp-construction} becomes
\[
 k'-W=\frac{t^2}{12}+\frac{td}{3}-d^2-\frac t3+\eta_2(n)
 =\frac{t^2}{9}-\left(d-\frac t6\right)^2-\frac t3+\eta_2(n).
\]
Choosing $d$ within $1/2$ of $t/6$ so that both part sizes are integers gives $k'-W=t^2/9+O(t+1)$. Since $t\sim2k/n$, this yields $A_3=4/9$.
\end{remark}

\section{Proof of Theorem~\ref{thm:uniform}}\label{sec:fractional}

We prove the uniform fractional bound in Theorem~\ref{thm:uniform}. The constant $c_r>0$ in that theorem will be chosen sufficiently small, depending only on $r$. Let $e(G)=t_{r-1}(n)+k$, where $0\le k\le c_rn^2$, and let $D=e(G)-\nu_r^*(G)-\frac{r-2}{2(r-1)}n^2$ denote the deficit. Since $k-\nu_r^*(G)=D+\eta_{r-1}(n)$ and the rounding term is bounded, it suffices to bound $D$ when $D>0$ and $n$ is sufficiently large. Starting with a minimum fractional cover, we delete its edges of weight one and find $r-1$ large independent sets in the remaining graph. We adjust the cover and enlarge the parts so that Lemma~\ref{lem:model} bounds $D$ in terms of the number of vertices left outside them. We bound this number by counting the edges of the remaining graph and comparing the result with $t_{r-1}(n)+k$.

\subsection{Finding large independent sets}

We use the following stability theorem of F\"uredi~\cite[Corollary~2]{Furedi} on the zero-weight edges of a minimum cover.

\begin{lemma}\label{lem:core}
Every $K_r$-free graph $F$ on $n$ vertices admits a complete $(r-1)$-partite graph $K$ on $V(F)$ such that
\[
 |E(F)\mathbin{\triangle}E(K)|\le3\bigl(t_{r-1}(n)-e(F)\bigr).
\]
\end{lemma}

The deficit $D$ is unchanged when we delete from $G$ the edges given weight one by a minimum fractional cover $w$. Let $G_0$ denote the resulting graph. Restricting $w$ to $G_0$ reduces its cost by one per deleted edge, and any cover of $G_0$ extends to $G$ by giving each deleted edge weight one. Thus the restriction is optimal and $e(G_0)-\nu_r^*(G_0)=e(G)-\nu_r^*(G)$. Let $k_0=e(G_0)-t_{r-1}(n)$. Then $0<k_0\le k$ and $\nu_r^*(G_0)=k_0-\eta_{r-1}(n)-D<k_0$.

The zero-weight edges of $G_0$ form a $K_r$-free spanning graph $F$. By complementary slackness, every other edge has load one in an optimal fractional packing, whose total edge load is $\binom r2\nu_r^*(G_0)$, so
\[
 e(F)\ge e(G_0)-\binom r2\nu_r^*(G_0)
 \ge t_{r-1}(n)-\left(\binom r2-1\right)k_0.
\]
By Lemma~\ref{lem:core}, $F$ differs in $O_r(k+1)$ edges from a complete $(r-1)$-partite graph $K$ with parts $V_1,\ldots,V_{r-1}$. Hence $e(K)\ge t_{r-1}(n)-O_r(k+1)$ and
\[
 \sum_i\left(|V_i|-\frac n{r-1}\right)^2
 =\frac{r-2}{r-1}n^2-2e(K)=O_r(k+1).
\]
For sufficiently small $c_r$, every $V_i$ has size at least $3n/[4(r-1)]$.

Let $S$ consist of the vertices missing at least $n/[8(r-1)^2]$ edges of $F$ to other parts of $K$, and let $A_i=V_i\setminus S$. Each missing edge in $E(K)\setminus E(F)$ contributes twice to the missing degrees. Thus
\begin{equation}\label{eq:initial-core}
 s=|S|=O_r((k+1)/n).
\end{equation}
Write $b_i=|A_i|$ and $b_{\min}=\min_i b_i$. Since $k\le c_rn^2$, equation~\eqref{eq:initial-core} gives $s=O_r(c_rn+1/n)$. Decreasing $c_r$ if necessary therefore gives $b_{\min}\ge n/[2(r-1)]$ and $s/b_{\min}\le1/4$.

Each $A_i$ is independent in $G_0$. Otherwise, extend an edge $uv\in E(G_0[A_i])$ to a $K_r$ by choosing one vertex from each other $A_j$, with all other edges in $F$. At each step the vertices already chosen exclude fewer than $(r-1)n/[8(r-1)^2]<b_j$ candidates. This clique would have total weight $w_{uv}<1$, contradicting the cover constraint. Also, $s\ge1$, since $G_0$ has more than $t_{r-1}(n)$ edges.

\subsection{A bound for the cover}

The following lemma allows the weights from $S$ to each part to vary with the vertex of $S$. It bounds the sum of $1-w_{uv}$ over all unordered vertex pairs and retains a sum of squares that controls the part sizes.

\begin{lemma}\label{lem:model}
Partition an $n$-vertex set into $S,A_1,\ldots,A_{r-1}$, where $s=|S|\ge1$
and $b_i=|A_i|>0$. Let $w$ assign a weight in $[0,1]$ to each pair so
that every $r$-set has total weight at least one. Suppose $w$ is one
inside each $A_i$ and zero between distinct parts $A_i,A_j$.
For $v\in S$, suppose every pair from $v$ to $A_i$ has weight $a_i(v)$,
where $\sum_i a_i(v)=1$. Let $\mu_i=s^{-1}\sum_{v\in S}a_i(v)$.
Then
\begin{equation}\label{eq:model-bound}
 \sum_{\{u,v\}}(1-w_{uv})-\frac{r-2}{2(r-1)}n^2\le \frac{A_r}{(r-1)^2}s^2
       -\frac12\sum_{i=1}^{r-1}(b_i-n/(r-1)+s\mu_i)^2 .
\end{equation}
\end{lemma}

\begin{proof}
For this calculation let $w_{vv}=1$ for $v\in S$, and let $z=(2s^2)^{-1}\sum_{u,v\in S}w_{uv}$, summing over ordered pairs. Thus $2z$ is the mean weight of two independent uniform choices from $S$. The unordered pairs inside $S$ contribute $s^2(1/2-z)$ to $\sum_{\{u,v\}}(1-w_{uv})$. Adding the other pairs and completing the square using $\sum_i b_i=n-s$ and $\sum_i\mu_i=1$ gives
\begin{align*}
 \sum_{\{u,v\}}(1-w_{uv})-\frac{r-2}{2(r-1)}n^2
 &=\frac{n^2}{2(r-1)}-zs^2-\frac12\sum_i b_i^2-s\sum_i b_i\mu_i\\
 &=s^2\left(-z+\frac12\sum_i\mu_i^2\right)
   -\frac12\sum_i\left(b_i-\frac n{r-1}+s\mu_i\right)^2.
\end{align*}
It remains to bound $-z+\tfrac12\sum_i\mu_i^2$ by $A_r/(r-1)^2$. We compare the means $\mu_i$ with the weights $x_i$ in~\eqref{eq:optimal-weights}.

Fix $2\le j\le r$ and a set $I$ of $r-j$ part indices. Choose one vertex from each indexed part and sample $j$ vertices independently and uniformly from $S$ with replacement. The total weight is at least one even if $s<j$, since choosing distinct vertices gives a covered $r$-set, and a repeated vertex gives a pair of weight one. Using the pair mean $2z$ and the part means $\mu_i$, we obtain
\begin{equation}\label{eq:model-prefix-source}
 j(j-1)z+j\sum_{i\in I}\mu_i\ge1.
\end{equation}
In particular, $z\ge1/[r(r-1)]$ by taking $j=r$.

Relabel the parts so that $\mu_1\le\cdots\le\mu_{r-1}$. For $1\le\ell\le r-2$, take $j=r-\ell$ and $I=\{1,\ldots,\ell\}$ in~\eqref{eq:model-prefix-source}. Comparing with~\eqref{eq:prefix-weights} yields
\[
 \sum_{i=1}^{\ell}(\mu_i-x_i)
 \ge -(r-\ell-1)\left(z-\frac1{r(r-1)}\right).
\]
Both sequences have sum one. Summation by parts and the monotonicity of $\mu_i+x_i$ give
\begin{align*}
 \sum_i(\mu_i^2-x_i^2)
 &=-\sum_{\ell=1}^{r-2}\left(\sum_{i=1}^{\ell}(\mu_i-x_i)\right)
       (\mu_{\ell+1}+x_{\ell+1}-\mu_\ell-x_\ell)\\
 &\le\left(z-\frac1{r(r-1)}\right)
       \sum_{\ell=1}^{r-2}(r-\ell-1)(\mu_{\ell+1}+x_{\ell+1}-\mu_\ell-x_\ell).
\end{align*}
The last sum telescopes to $2-(r-1)(\mu_1+x_1)\le2$. Hence
\[
 -z+\frac12\sum_i\mu_i^2
 \le-\frac1{r(r-1)}+\frac12\sum_i x_i^2
 =\frac{A_r}{(r-1)^2}.
\]
Substituting this inequality into the completed-square identity proves~\eqref{eq:model-bound}.
\end{proof}

To apply Lemma~\ref{lem:model}, extend the minimum cover $w$ of $G_0$ to all unordered pairs of distinct vertices by giving every nonedge weight one. Every $r$-set is covered, since it spans a $K_r$ or contains a nonedge. The cost on $G_0$ remains $\sum_{uv\in E(G_0)}w_{uv}$. Since the nonedges contribute zero to $\sum_{\{u,v\}}(1-w_{uv})$, this sum equals $e(G_0)-\nu_r^*(G_0)$. Hence $D=\sum_{\{u,v\}}(1-w_{uv})-\frac{r-2}{2(r-1)}n^2$, and lowering the total pair weight increases this deficit by the same amount.

The average over permutations within each part gives a feasible cover of the same total weight. In the original extended cover, let $L_{ij}=\sum_{x\in A_i,\,y\in A_j}w_{xy}$ for $i\ne j$, $L=\sum_{i<j}L_{ij}$, and $\alpha_i(v)=b_i^{-1}\sum_{x\in A_i}w_{vx}$ for $v\in S$. For $U\subseteq S$ and a set $I$ of part indices with $|U|+|I|=r$, the averaged constraint is
\begin{equation}\label{eq:averaged-cover}
 \sum_{\{u,v\}\subseteq U}w_{uv}
 +\sum_{v\in U}\sum_{i\in I}\alpha_i(v)
 +\sum_{\substack{i<j\\i,j\in I}}\frac{L_{ij}}{b_ib_j}\ge1.
\end{equation}
The last term is lost when the weights between parts are set to zero. For each $v\in S$ and $i<j$, we compensate by adding half the mean weight between $A_i$ and $A_j$ to each of $\alpha_i(v)$ and $\alpha_j(v)$. Let
\[
 \delta_v=\sum_i\alpha_i(v)+\sum_{i<j}\frac{L_{ij}}{b_ib_j}-1.
\]
We then normalize the resulting weights by setting
\[
 a_i(v)=\frac{\alpha_i(v)+\tfrac12\sum_{j\ne i}L_{ij}/(b_ib_j)}{1+\delta_v}.
\]
By~\eqref{eq:averaged-cover} with $U=\{v\}$, we have $\delta_v\ge0$. The numerators are nonnegative and sum to $1+\delta_v$, so $\sum_i a_i(v)=1$.

The normalized weights are one inside each $A_i$, zero between parts, and $a_i(v)$ on pairs from $v$ to $A_i$. Each pair $uv$ in $S$ receives weight $\min\{1,w_{uv}+\delta_u+\delta_v\}$. Any $r$-set meeting a part twice is covered. Otherwise it meets distinct parts indexed by $I$ and a nonempty set $U\subseteq S$, since there are only $r-1$ parts. Before division, each $v\in U$ gains
$\tfrac12\sum_{i\in I}\sum_{j\ne i}L_{ij}/(b_ib_j)\ge\sum_{i<j,\ i,j\in I}L_{ij}/(b_ib_j)$
on its pairs to these parts, at least the weight lost between the parts. The division then removes at most $\delta_v$ from these pairs. If $|U|\ge2$, the increase within $U$ is $(|U|-1)\sum_{v\in U}\delta_v$ unless a pair reaches weight one. In both cases the set remains covered. If $|U|=1$, all external parts occur and the new total weight is $\sum_i a_i(v)=1$.

The total weight between the parts drops by $L$. On pairs from $S$ to the parts, the increases before division cost at most
$\frac{s}{2}\sum_{i<j}L_{ij}(b_i^{-1}+b_j^{-1})\le sL/b_{\min}$.
Division decreases the sum of the entries for each $v$ by $\delta_v$, saving at least $b_{\min}\sum_v\delta_v$ in pair weight. The increase inside $S$ costs at most $(s-1)\sum_v\delta_v$. Comparing these savings with the upper bound on the new deficit in Lemma~\ref{lem:model} gives
\[
 D+(1-s/b_{\min})L+(b_{\min}-s+1)\sum_v\delta_v
 \le\frac{A_r}{(r-1)^2}s^2-\frac12\sum_i\left(b_i-\frac n{r-1}+\sum_{v\in S}a_i(v)\right)^2.
\]
Since $D>0$, $0\le\sum_v a_i(v)\le s$, $b_{\min}\ge n/[2(r-1)]$, and $s/b_{\min}\le1/4$, it follows that
\begin{equation}\label{eq:normalization-bounds}
\begin{gathered}
 L=O_r(s^2),\\
 \sum_v\delta_v=O_r(s^2/n),\\
 b_i=n/(r-1)+O_r(s).
\end{gathered}
\end{equation}
For later degree estimates in $G_0$, the original means satisfy $\alpha_i(v)\le (1+\delta_v)a_i(v)\le a_i(v)+\delta_v$. This inequality and the definition of $\delta_v$ give
\begin{equation}\label{eq:original-arms}
\begin{gathered}
 1-\alpha_i(v)\ge1-a_i(v)-\delta_v,\\
 \sum_i(1-\alpha_i(v))=r-2-\delta_v+\sum_{i<j}\frac{L_{ij}}{b_ib_j}.
\end{gathered}
\end{equation}

\subsection{Enlarging the parts}

A vertex of $S$ can be assigned to a part at small cost when one of its normalized weights is close to one. For sufficiently small $c_r$ we have $\eta=\sqrt{s/n}<1/2$. Let $I=\{v\in S:\max_i a_i(v)\le1-\eta\}$ and $u=|I|$. The degrees in $G_0$ bound the number $u$ of vertices that will remain outside the enlarged parts.

\begin{claim}\label{clm:effective-core}
\begin{equation}\label{eq:effective-core}
 u\le\frac{(r-1)k_0}{n}+O_r(s\eta).
\end{equation}
\end{claim}

\begin{proof}
We show that the vertices of $I$ have few nonneighbors outside $S$ in total. For $v\in I$, let $H_v=\sum_i|A_i\setminus N_{G_0}(v)|$. We use the original cover of $G_0$, extended by weight one on nonedges. By~\eqref{eq:normalization-bounds}, $\sum_{j<\ell}L_{j\ell}/(b_jb_\ell)\le L/b_{\min}^2=O_r(s^2/n^2)=O_r(\eta^4)$. Since $a_i(v)\le1-\eta$, equation~\eqref{eq:original-arms} therefore gives
\[
 \sum_{j\ne i}\alpha_j(v)
 =1+\delta_v-\alpha_i(v)-\sum_{j<\ell}\frac{L_{j\ell}}{b_jb_\ell}
 \ge\eta-O_r(\eta^4)\ge\eta/2
\]
for sufficiently small $c_r$.

The mean total weight on $v$ and one independent uniform choice from each part is $1+\delta_v$, and every such $r$-set has weight at least one. For a fixed $x\in A_i\setminus N_{G_0}(v)$, the pair $vx$ has weight one. Averaging over the other parts, the mean weight in excess of one is therefore at least $\sum_{j\ne i}\alpha_j(v)$. All remaining choices have nonnegative excess, so
\[
 \delta_v\ge\frac{|A_i\setminus N_{G_0}(v)|}{b_i}\sum_{j\ne i}\alpha_j(v)
 \ge\frac{\eta|A_i\setminus N_{G_0}(v)|}{2b_i}.
\]
Summing over $i$ gives $H_v\le2\delta_v\sum_i b_i/\eta\le2n\delta_v/\eta$. Thus
\[
 \sum_{v\in I}H_v\le\frac{2n}{\eta}\sum_{v\in S}\delta_v
 =O_r(s^2/\eta)=O_r(sn\eta).
\]

For $v\in S\setminus I$, its degree outside $S$ is at least the total of $1-w_{vx}$ over the external vertices $x$. By~\eqref{eq:normalization-bounds} and~\eqref{eq:original-arms},
\[
 |N_{G_0}(v)\setminus S|\ge\sum_i b_i(1-\alpha_i(v))
 \ge\frac{r-2}{r-1}n-\frac n{r-1}\delta_v-O_r(s).
\]
For $v\in I$ the degree outside $S$ is $n-s-H_v$. Between distinct parts, every missing edge contributes one to $L$. Since their sizes differ from $(n-s)/(r-1)$ by $O_r(s)$, the parts span at least $\sum_{i<j}b_ib_j-L=\frac{r-2}{2(r-1)}(n-s)^2-O_r(s^2)$ edges. Adding the edges from $S$ and using $n\sum_v\delta_v=O_r(s^2)$ gives
\begin{align*}
 e(G_0)&\ge\frac{r-2}{2(r-1)}(n-s)^2
       +\frac{r-2}{r-1}n(s-u)+u(n-s)-\sum_{v\in I}H_v-O_r(s^2)\\
 &\ge\frac{r-2}{2(r-1)}n^2+\frac{un}{r-1}-O_r(sn\eta).
\end{align*}
Here $s^2\le sn\eta$. Comparing with $e(G_0)=\frac{r-2}{2(r-1)}n^2+k_0-\eta_{r-1}(n)$ proves the claim.
\end{proof}

We now enlarge the parts in the normalized cover, leaving only $I$ outside them. The enlarged parts need not be independent in $G_0$. Each $v\in S\setminus I$ has a unique index $i$ with $a_i(v)>1-\eta$, since the weights sum to one and $\eta<1/2$. Assign $v$ to $A_i$ and let $\gamma_v=1-a_i(v)$. Give pairs inside an enlarged part weight one, pairs between enlarged parts weight zero, and pairs from the enlarged $A_i$ to $z\in I$ weight $a_i(z)$. Keep the weights on $I$.

An $r$-set meeting an enlarged part twice contains a pair of weight one. Otherwise replace each newly assigned vertex by any vertex of the original part to which it was assigned. The distinct replacements give an $r$-set in the normalized cover whose pair weights equal those of the original $r$-set in the new cover. Thus the new weights form a cover.

For $v$ assigned to $A_i$, its pairs to the original parts change in total weight by
\[
 b_i\gamma_v-\sum_{j\ne i}b_ja_j(v)
 =\sum_{j\ne i}(b_i-b_j)a_j(v)=O_r(s\gamma_v).
\]
To control the change inside $S$, take $v$, another vertex $z\in S$, and one vertex from each original part other than $A_i$. If $z\in I$, the pairs to these parts have total normalized weight $\gamma_v+1-a_i(z)$, so $vz$ has weight at least $a_i(z)-\gamma_v$. If $z$ is also assigned to $A_i$, the same constraint gives weight at least $1-\gamma_v-\gamma_z$ on $vz$. A pair whose endpoints are assigned to distinct parts receives weight zero. The total increase inside $S$ is therefore at most $(s-1)\sum_{v\in S\setminus I}\gamma_v$. Together with the change on pairs to the original parts, this gives a total increase of at most $O_r(\eta s^2)$.

The normalized cover had deficit at least $D$, so the new deficit is at least $D-O_r(\eta s^2)$. If $u>0$, Lemma~\ref{lem:model} bounds it above by $A_ru^2/(r-1)^2$. If $u=0$, the pairs of weight zero form a complete $(r-1)$-partite graph with at most $\frac{r-2}{2(r-1)}n^2$ edges, so the deficit is nonpositive. Thus
\begin{equation}\label{eq:fractional-close}
 D\le\frac{A_r}{(r-1)^2}u^2+O_r(\eta s^2).
\end{equation}

\smallskip
\begin{proof}[Proof of Theorem~\ref{thm:uniform}]
For $D\le0$ the bound follows from $0\le\eta_{r-1}(n)\le(r-1)/8$. If $D>0$ and $k<n$, equation~\eqref{eq:initial-core} gives $s=O_r(1)$, and the normalized cover gives $D\le A_rs^2/(r-1)^2=O_r(1)$.

It remains to consider $D>0$ and $k\ge n$. Then~\eqref{eq:initial-core} gives $s=O_r(k/n)$ and $\eta=O_r(\sqrt{k}/n)$. Since $k_0\le k$, equation~\eqref{eq:effective-core} gives $u\le(r-1)k/n+O_r(k^{3/2}/n^2)$. Substituting in~\eqref{eq:fractional-close} yields
\[
 D\le\frac{A_r}{(r-1)^2}
 \left(\frac{(r-1)k}{n}+O_r\!\left(\frac{k^{3/2}}{n^2}\right)\right)^2
 +O_r\!\left(\frac{k^{5/2}}{n^3}\right)
 =A_r\frac{k^2}{n^2}+O_r\!\left(\frac{k^{5/2}}{n^3}\right).
\]
Here $k\le c_rn^2$ gives $k^3/n^4=O_r(k^{5/2}/n^3)$. Adding the bounded term $\eta_{r-1}(n)$ proves the lower bound on $\nu_r^*(G)$ for sufficiently large $n$. Increasing $B_r$ covers the finitely many smaller orders.
\end{proof}

\smallskip
\begin{proof}[Proof of Theorem~\ref{thm:second-order} (Fractional Case)]
For $n\ll k\ll n^2$, the error terms $k^{5/2}/n^3$ and $1$ are $o(k^2/n^2)$. The result follows from Theorem~\ref{thm:uniform} and~\eqref{eq:sharp-upper}.
\end{proof}

\section{Proof of Theorem~\ref{thm:second-order} (Integer Case)}\label{sec:integer}

For the integer lower bound, let $G$ be an $n$-vertex graph with $e(G)=t_{r-1}(n)+k$, where $n\ll k\ll n^2$. We may assume $\nu_r(G)<k$. The theorem of Haxell and R\"odl~\cite{HR} gives $\nu_r^*(J)-\nu_r(J)=o_r(N^2)$ uniformly over graphs $J$ on $N$ vertices. See also~\cite[Theorem~1.1]{Yuster}. The error $o_r(n^2)$ does not suffice here. We obtain the required error $o(k^2/n^2)$ by reducing to a graph on $O_r(k/n)$ vertices.

We first modify $G$ to obtain an $n$-vertex graph $G_1$ whose excess $k_1$ satisfies $k-\nu_r(G)\le k_1-\nu_r(G_1)$. Section~\ref{sec:fractional} bounds $k_1-\nu_r^*(G_1)$. We bound the remaining difference $\nu_r^*(G_1)-\nu_r(G_1)$ by that in the smaller graph.

\subsection{Removing edges inside the parts}

The edges not used by a maximum $K_r$-packing in $G$ form a $K_r$-free spanning graph $F$ with at least $t_{r-1}(n)-(\binom r2-1)k$ edges. F\"uredi's stability theorem~\cite[Corollary~2]{Furedi} gives a complete $(r-1)$-partite graph differing from $F$ in $O_r(k)$ edges. Since $e(G)-e(F)=O_r(k)$, the corresponding partition $A_1,\ldots,A_{r-1}$ also has only $O_r(k)$ internal edges and missing edges between parts in $G$. In particular,
\[
 \frac12\sum_i\left(|A_i|-\frac n{r-1}\right)^2
 =\frac{r-2}{2(r-1)}n^2-\sum_{i<j}|A_i||A_j|=O_r(k),
\]
so $|A_i|=n/(r-1)+o(n)$.

Fix $\varepsilon=\varepsilon(r)>0$ with $r(2r-1)\varepsilon<1/[2(r-1)]$. Let $S$ consist of the vertices for which the number of neighbors in their own part plus the number of nonneighbors in the other parts is at least $\varepsilon n$. Each internal edge and each missing edge between parts contributes twice, so $s=|S|=O_r(k/n)$. Let $B_i=A_i\setminus S$ and $b_i=|B_i|$. These external parts have sizes $b_i\sim n/(r-1)$. Each vertex of $B_i$ has fewer than $\varepsilon n$ internal neighbors and fewer than $\varepsilon n$ nonneighbors in the other parts combined. Let $L$ be the set of edges inside the $B_i$ and let $G_0=G-L$.

\begin{lemma}\label{lem:absorb}
For the partition above and sufficiently large $n$, $\nu_r(G)\ge|L|+\nu_r(G_0)$.
\end{lemma}

\begin{proof}
Fix a packing $\mathcal P$ in $G_0$. Every clique of $\mathcal P$ meets $S$, since the $B_i$ are independent. The cliques containing a given external vertex use distinct edges from it to $S$. Thus there are at most $s$ such cliques, using at most $(r-2)s$ of its edges to other external parts.

Extend each $uv\in L$ in turn by choosing one vertex from each of the other external parts, avoiding edges used by $\mathcal P$ or any earlier extension. Allow each vertex to be chosen as an auxiliary vertex at most $\lceil\varepsilon n\rceil$ times. Since it is also an endpoint of fewer than $\varepsilon n$ edges of $L$, the cliques added so far use at most $(r-1)(2\varepsilon n+1)$ of its edges. When we choose a vertex in $B_i$, the at most $r$ vertices already in the current clique exclude at most
\[
 r\bigl((2r-1)\varepsilon n+(r-2)s+r-1\bigr)
\]
candidates through missing or used edges. At most $(r-2)|L|/\lceil\varepsilon n\rceil=O_r(k/n)$ further vertices have reached their limit on choices. Since $s=O_r(k/n)=o(n)$ and $r(2r-1)\varepsilon<1/[2(r-1)]$, fewer than $b_i\sim n/(r-1)$ candidates are excluded for large $n$. Each extension therefore exists. Its only internal edge is the edge being processed, so all $|L|$ new cliques are edge-disjoint and avoid $\mathcal P$.
\end{proof}

Let $k_0=k-|L|$. Then $e(G_0)=t_{r-1}(n)+k_0$ and Lemma~\ref{lem:absorb} gives $k-\nu_r(G)\le k_0-\nu_r(G_0)$. In particular, $\nu_r(G_0)<k_0$ and $k_0>0$. Since $G_0-S$ is $(r-1)$-partite, this positive excess implies $s\ge1$.

\subsection{Completing missing edges}

We next add the missing edges between the $B_i$ and from each $v\in S$ to all but one of the parts. To bound their number, we compare $k_0$ with a packing in which each clique uses one vertex of $S$. The following lemma supplies the required external cliques.

\begin{lemma}\label{lem:transversal}
Start with a complete $(r-1)$-partite graph with part sizes $p_1,\ldots,p_{r-1}$. Delete at most $M$ edges and then the edges of a graph of maximum degree at most $D$. The remaining graph has at least $\min_i p_i-O_r(D+\sqrt M+1)$ vertex-disjoint copies of $K_{r-1}$.
\end{lemma}

\begin{proof}
Trim the parts to size $\min_i p_i$ and choose a maximal vertex-disjoint $K_{r-1}$-packing. Let $m$ vertices remain in each part. If $m>0$, every choice of one remaining vertex from each part contains a missing edge, by maximality. Each missing edge belongs to $m^{r-3}$ of the $m^{r-1}$ choices. The degree bound accounts for at most $(r-1)mD/2$ missing edges in addition to the first $M$. Hence
\[
 m^{r-1}\le\left(M+\frac{(r-1)mD}{2}\right)m^{r-3}.
\]
Dividing by $m^{r-3}$ gives $m^2\le M+(r-1)mD/2$. It follows that $m\le(r-1)D/2+\sqrt M$. This also holds when $m=0$ and proves the lemma.
\end{proof}

In $G_0$, let $M$ count the missing edges between the $B_i$. A copy of $K_r$ meeting $S$ only in $v$ needs a neighbor of $v$ in every part. We therefore assign each $v\in S$ to a part $B_{j(v)}$ in which it has the fewest neighbors. Let $m_v=|N_{G_0}(v)\cap B_{j(v)}|$ and let
\[
 T=\sum_{v\in S}\sum_{i\ne j(v)}\bigl(b_i-|N_{G_0}(v)\cap B_i|\bigr).
\]
Thus $T$ counts the missing edges from vertices of $S$ to parts other than their assigned parts.

\begin{claim}\label{clm:integer-defects}
With this notation,
\begin{equation}\label{eq:integer-defects}
\begin{gathered}
 M=O_r(s^2),\\
 T=O_r(s^2),\\
 b_i=n/(r-1)+O_r(s).
\end{gathered}
\end{equation}
\end{claim}

\begin{proof}
We build a packing one vertex of $S$ at a time. At each $v$, the smallest part of its external neighborhood has size $m_v$. The copies of $K_{r-1}$ chosen for each earlier vertex were vertex-disjoint, so the previously used external edges have maximum degree at most $(r-2)s$. Lemma~\ref{lem:transversal} applies to the neighborhood after we delete these edges. Extending the resulting $m_v-O_r(s+\sqrt M+1)$ copies of $K_{r-1}$ by $v$ gives
\begin{equation}\label{eq:integer-baseline}
 \nu_r(G_0)\ge\sum_{v\in S}m_v-O_r(s^2+s\sqrt M).
\end{equation}

The assignment of $S$ to the parts also gives an exact edge count. Let $t_i=|\{v\in S:j(v)=i\}|$ and consider the complete $(r-1)$-partite graph whose $i$th part consists of $B_i$ and the vertices assigned to it. It has $\frac{r-2}{2(r-1)}n^2-\frac12\sum_i(b_i+t_i-n/(r-1))^2$ edges. To obtain $G_0$, delete the $M+T$ missing edges with at least one external endpoint and add the $\sum_v m_v$ edges to assigned parts. Within $S$, replace the $(s^2-\sum_i t_i^2)/2$ edges between parts by the edges of $G_0[S]$. Thus
\begin{equation}\label{eq:integer-squares}
 k_0-\sum_{v\in S}m_v
 =e(G_0[S])-M-T-\frac{s^2-\sum_i t_i^2}{2}
 -\frac12\sum_i\left(b_i+t_i-\frac n{r-1}\right)^2+\eta_{r-1}(n).
\end{equation}
Since $\nu_r(G_0)<k_0$, equation~\eqref{eq:integer-baseline} bounds $\sum_vm_v-k_0$ from above by $O_r(s^2+s\sqrt M)$. Using $e(G_0[S])\le\binom s2$, $\sum_i t_i^2\le s^2$, and $\eta_{r-1}(n)=O_r(1)$ in~\eqref{eq:integer-squares}, we obtain
\[
 M+T+\frac12\sum_i\left(b_i+t_i-\frac n{r-1}\right)^2\le C(s^2+s\sqrt M)
\]
for a constant $C=C(r)$. The inequality $Cs\sqrt M\le M/2+C^2s^2/2$ gives $M,T=O_r(s^2)$ and $|b_i+t_i-n/(r-1)|=O_r(s)$. Since $0\le t_i\le s$, the estimate for $b_i$ follows.
\end{proof}

Add the $M$ missing edges between the $B_i$ and the $T$ missing edges from vertices of $S$ to parts other than their assigned parts. If $b_i-m_v<2(r-1)^2s$ for some $i\ne j(v)$, also join $v$ to all of $B_{j(v)}$. By~\eqref{eq:integer-defects}, this adds at most $b_{j(v)}-m_v\le|b_{j(v)}-b_i|+2(r-1)^2s=O_r(s)$ edges per such vertex. The resulting graph $G_1$ has $a=O_r(s^2)$ more edges than $G_0$, and $G_1[S]=G_0[S]$.

In $G_1$, let $U$ consist of the vertices of $S$ adjacent to every external vertex, and let $X_j=\{v\in S\setminus U:j(v)=j\}$. Each $v\in X_j$ has $m_v$ neighbors in $B_j$, is joined to every other part, and satisfies
\begin{equation}\label{eq:slack}
 m_v\le\min_{i\ne j}b_i-2(r-1)^2s.
\end{equation}
This leaves enough vertices in every other part to extend all the edges from $v$ to $B_j$ to edge-disjoint cliques.

The excess of $G_1$ is $k_1=k_0+a$. Deleting the $a$ added edges destroys at most $a$ cliques in any packing, so $\nu_r(G_1)\le\nu_r(G_0)+a$. Together with Lemma~\ref{lem:absorb} and $a=O_r(s^2)=O_r(k^2/n^2)$, this gives
\begin{equation}\label{eq:integer-deficit-transfer}
\begin{gathered}
 k-\nu_r(G)\le k_1-\nu_r(G_1),\\
 0<k_1\le k+O_r(k^2/n^2)=k+o(k).
\end{gathered}
\end{equation}

\subsection{Reduction to a smaller graph}

The external parts of $G_1$ differ in size by only $O_r(s)$, so removing equally many vertices from each can leave $O_r(s)$ vertices in every part. We keep $S$ and delete the remaining edges from each $X_j$ to $B_j$. The next lemma bounds $\nu_r^*(G_1)-\nu_r(G_1)$ by the corresponding difference in the smaller graph.

\begin{lemma}\label{prop:compression}
Let $G_1$ have independent parts $B_1,\ldots,B_{r-1}$ spanning a complete $(r-1)$-partite graph and remaining vertices $S=U\cup X_1\cup\cdots\cup X_{r-1}$, where the sets in this union are disjoint and $s=|S|\ge1$. Each vertex of $U$ is adjacent to every external vertex. Each $v\in X_j$ has $m_v$ neighbors in $B_j$, is joined to every other part, and satisfies~\eqref{eq:slack}. Let $b_i=|B_i|$ and $b=\min_i b_i$. If $b\ge(r+3)s$, there is a graph $J$ on $(4r-3)s+\sum_i(b_i-b)$ vertices such that
\[
 \nu_r^*(G_1)=\sum_{v\in S\setminus U}m_v+|U|(b-4s)+\nu_r^*(J).
\]
Moreover, the integer packing number satisfies
\[
 \nu_r(G_1)\ge\sum_{v\in S\setminus U}m_v+|U|(b-4s)+\nu_r(J).
\]
In particular, $\nu_r^*(G_1)-\nu_r(G_1)\le\nu_r^*(J)-\nu_r(J)$.
\end{lemma}

\begin{proof}
The graph $J$ retains $H=G_1[S]$ and $4s+b_i-b$ vertices from each $B_i$, with all edges between the retained parts. Between $S$ and these parts, its edges are precisely those from $U$ to every part and from each $X_j$ to the parts other than part $j$. Thus $J$ is a subgraph of $G_1$ of the stated order.

For the fractional identity, we first weight the cliques of $J$ meeting $S$ at least twice, then complete the packing in each graph with cliques meeting $S$ once. Every copy $K$ of $K_r$ in $J$ satisfies
\begin{equation}\label{eq:allowed-types}
 V(K)\cap X_j\ne\varnothing\Longrightarrow V(K)\cap B_j=\varnothing.
\end{equation}
With $V(K)\cap S$ fixed, any choice of one vertex from each external part met by $K$ gives a clique in either graph. A nonnegative weight $w_K$ on each $K$ in $J$ meeting $S$ at least twice can therefore be distributed uniformly over these choices in either graph. The only constraints we impose on the $w_K$ are the edge capacities in $H$. For every $e\in E(H)$, we require
\begin{equation}\label{eq:core-capacity}
 \sum_{K:\,e\in E(K)}w_K\le1.
\end{equation}
For $v\in S$, let $l_{vi}$ be the sum of $w_K$ over cliques $K$ containing $v$ and meeting $B_i$. Each such clique uses an edge of $H$ incident to $v$. Summing these edge capacities gives
\begin{equation}\label{eq:core-loads}
 l_{vi}\le\sum_{K:\,v\in V(K)}w_K\le s-1.
\end{equation}

At a vertex $v\in U$ in $G_1$, the remaining total capacity to part $i$ is $b_i-l_{vi}$. We distribute total weight $\min_i(b_i-l_{vi})$ uniformly over cliques meeting $S$ only in $v$. At a vertex $v\in X_j$ we distribute total weight $m_v$ in the same way, using its actual neighbors in $B_j$ and adding nothing if $m_v=0$. The construction in $J$ uses $4s+b_i-b$ in place of $b_i$ at $v\in U$. No clique in $J$ meets $S$ only in a vertex of $X_j$, since it would need a neighbor in part $j$. The added weights are nonnegative by~\eqref{eq:core-loads}.

These additions leave the loads on $H$ unchanged, so they satisfy~\eqref{eq:core-capacity}. At each $v\in U$, the added weight is at most the remaining capacity to each part and is spread uniformly. For $v\in X_j$ in $G_1$, the initial load on edges to $B_j$ is zero by~\eqref{eq:allowed-types}. Distributing the added weight $m_v$ uniformly over the $m_v$ such edges gives load one on each edge. Edges to $B_i$ for $i\ne j$ have load $(m_v+l_{vi})/b_i\le1$ by~\eqref{eq:slack} and~\eqref{eq:core-loads}. For $v\in X_j$ in $J$, all edges from $v$ to the parts have load at most $(s-1)/(4s)$.

For an edge $xy$ between external parts, fix $v\in S$. If $v\in X_j$, exchange $x$ and $y$ if necessary so that the part containing $y$ is not $B_j$. If $vx$ is a nonedge, no clique contains all three vertices. Otherwise its load is at most one by the preceding paragraph. Among cliques containing $v$ and $x$ and using the part containing $y$, both the initial and added weights are distributed uniformly over that whole part, independently of the other external choices. The total weight of cliques containing $v,x,y$ is therefore at most the load on $vx$ divided by the size of that part, and hence at most $1/(4s)$. Every clique meets $S$, so summing over $v$ bounds the load on $xy$ by $s/(4s)=1/4$.

These weights give a feasible fractional packing in each graph. We next show that each graph has an optimum of this form. In an optimal fractional packing of $G_1$, the cliques using an edge from some $X_j$ to $B_j$ have total weight at most $\sum_{v\in S\setminus U}m_v$. Once we remove them, each remaining clique meeting $S$ at least twice satisfies~\eqref{eq:allowed-types}. Its weight can therefore be transferred to a clique $K$ of $J$ with the same vertices in $S$ and the same external parts. We add all weights assigned to the same $K$, which preserves the loads in $H$ and the sums $l_{vi}$. The surviving cliques meeting $S$ once do so at vertices of $U$. For each $v\in U$, summing the edge capacities from $v$ to each part bounds their total weight at $v$ by $\min_i(b_i-l_{vi})$. Completion in $G_1$ thus restores at least the original packing weight. The same argument applies in $J$ with its smaller part sizes and no deletion. Since
\[
 \min_i(4s+b_i-b-l_{vi})=\min_i(b_i-l_{vi})-(b-4s),
\]
the total packing weight in $G_1$ equals that in $J$ plus $\sum_{v\in S\setminus U}m_v+|U|(b-4s)$ for every fixed choice of the $w_K$. Maximizing over the $w_K$ proves the fractional identity.

For the integer inequality, embed a maximum packing of $J$ in $G_1$. To add cliques through $U$, label the $b-4s$ vertices outside $J$ in each external part by $\mathbb Z_{b-4s}$ and number the parts $0,\ldots,r-2$. For each $0\le\ell<|U|$, the vertex tuples $(x,x+\ell,\ldots,x+(r-2)\ell)$ with $x\in\mathbb Z_{b-4s}$ and coordinates in the corresponding parts form a $K_{r-1}$-factor, since addition by $i\ell$ permutes the labels in part $i$. Distinct factors cannot share an edge between parts $i<j$, since that would give $(j-i)(\ell-\ell')\equiv0\pmod{b-4s}$, whereas $0<|(j-i)(\ell-\ell')|\le(r-2)(|U|-1)<b-4s$. Here $b-4s\ge(r-1)s$. Assign a different factor to each vertex of $U$ and extend its cliques by that vertex. This adds $(b-4s)|U|$ edge-disjoint cliques outside the embedded packing.

Every edge from $X_j$ to $B_j$ is still unused. For each $v\in X_j$ in turn, we extend all $m_v$ of these edges at $v$, choosing one vertex from every other part and using distinct vertices in each part for this fixed $v$. When we choose a vertex in $B_i$, the embedded packing excludes at most $s-1$ candidates at $v$, since every embedded clique through $v$ uses a distinct incident edge of $H$. Earlier choices for this $v$ exclude at most $m_v-1$ more. Each already chosen external vertex belongs to at most $s$ previously packed cliques, each using a distinct edge from it to $S$, so it has at most $s$ used edges to $B_i$. There are at most $r-2$ such vertices, leaving at least
\[
 b_i-(s-1)-(m_v-1)-(r-2)s=b_i-m_v-(r-1)s+2>0
\]
candidates by~\eqref{eq:slack}. This adds $\sum_{v\in S\setminus U}m_v$ cliques and proves the integer inequality.
\end{proof}

\smallskip
\begin{proof}[Proof of Theorem~\ref{thm:second-order} (Integer Case)]
The reductions above give $G_1$ satisfying~\eqref{eq:integer-deficit-transfer}, with $s=O_r(k/n)$ and $b_i-b=O_r(s)$ by~\eqref{eq:integer-defects}. Since $b_i\sim n/(r-1)$ and $s=o(n)$, Lemma~\ref{prop:compression} applies and gives a graph $J$ on $O_r(k/n)$ vertices. Add isolated vertices until its order is $\lceil Ck/n\rceil$, for a sufficiently large fixed $C=C(r)$. This changes neither packing number and ensures that the order tends to infinity, even if the original $J$ has bounded order. The Haxell--R\"odl theorem~\cite{HR} therefore gives, uniformly in $G$,
\[
 0\le\nu_r^*(G_1)-\nu_r(G_1)\le\nu_r^*(J)-\nu_r(J)
 =o_r\bigl(\lceil Ck/n\rceil^2\bigr)=o_r(k^2/n^2).
\]
The excess $k_1$ may be smaller than $n$, but $0<k_1\le k+O_r(k^2/n^2)=o(n^2)$. Hence the uniform bound of Theorem~\ref{thm:uniform} applies to $G_1$, and
\begin{align*}
 k-\nu_r(G)&\le k_1-\nu_r(G_1)
 =k_1-\nu_r^*(G_1)+\nu_r^*(G_1)-\nu_r(G_1)\\
 &\le A_r\frac{k_1^2}{n^2}+B_r\frac{k_1^{5/2}}{n^3}+B_r+o_r(k^2/n^2).
\end{align*}
Now $k_1/n\le(1+o(1))k/n$ and $k_1/n^2=o(1)$, so $k_1^{5/2}/n^3=(k_1/n)^2\sqrt{k_1/n^2}=o(k^2/n^2)$. Also $B_r=o(k^2/n^2)$ because $k/n\to\infty$. Thus $k-\nu_r(G)\le(A_r+o(1))k^2/n^2$ for every $G$ under consideration. The construction in~\eqref{eq:sharp-upper} gives the matching upper bound on $\phi_r(n,k)$ for every feasible edge count in this range.
\end{proof}

\section*{Acknowledgments}
Jialin He is supported in part by the National Natural Science Foundation of China grant~12601663 and the Science and Technology Commission of Shanghai Municipality grant~22DZ2229014. Fuhong Ma is supported by the National Natural Science Foundation of China grant~12401458 and the Natural Science Foundation of Shandong Province grant~ZR2024QA239. Jin Yan is supported by the National Natural Science Foundation of China grant~12571373 and the Natural Science Foundation of Shandong Province grant~ZR2025MS05.

\section*{Declaration on the Use of AI}
GPT-6 Astra assisted with discussions of the construction and its fractional cover. OpenAI Codex assisted with mathematical exploration, proof checking, literature searches, and preparation of the work on the second-order term. The tools also assisted with proofreading and language editing. The authors take full responsibility for all AI-assisted work and for the content of the manuscript.

\end{document}